\documentclass[12pt]{article}
\usepackage[utf8]{inputenc}
\usepackage[T1]{fontenc}

\usepackage{palatino}

\usepackage{amsfonts}

\usepackage{mdframed}

\usepackage{amssymb, amsmath}
\usepackage[english]{babel}
\usepackage{graphicx}
\graphicspath{ {./images/} }
\usepackage{subcaption}
\usepackage[export]{adjustbox}
\usepackage{wrapfig}
\usepackage{url}
\usepackage{setspace} 

\title{Thirty years after: insights on the cultural origins of Andrej N. Kolmogorov's 1954 invariant tori theorem from a short conversation with Vladimir I. Arnold  }
\author{Fascitiello Isabella}

\date{\today}

\begin{document}
\maketitle
\tableofcontents

\begin{abstract}
\noindent
Among the impressive contributions of Andrej N. Kolmogorov's to mathematics in the 20th century, his 1954 invariant tori theorem is still little understood from a historical point of view [Dumas 2014]. Vladimir I. Arnold – who entered Moscow State University that same year and would became a student of Kolmogorov – 30 years after asked him about the origin of a research program to which Arnold himself gave a big contribution. In this paper we compare and put in historical and biographical context two different account's by Arnold on this conversation. 

\end{abstract}
\noindent

\section{Introduction}

Four USSR mathematicians took part at the 1954 International Congress of Mathematicians at Amsterdam, the second one held after the end of the Second World War. Andrej N Kolmogorov, one of them, has been invited to close the Congress with plenary conference devoted to the theory of dynamical systems and classical mechanics. He presented a research program that was linked to contributions in this area in the 1930s – a period of intense scientific activity in countries where the nightmare of war and totalitarian terror was about to break down. In a short note written for the 1985 volume of his  Selected works (edited by his former student Vladimir Mikhailovich Tikhomirov), he wrote: 
\\

\leftskip=1cm
\noindent
My papers on classical mechanics appeared under the influence of von Neumann's papers on the spectral theory of dynamical systems and, particularly under the influence of the Bogolyubov-Krylov paper of 1937. I became extremely interested in the question of what ergodic sets (in the sense of Bogolyubov-Krylov) can exist in the dynamical systems of classical mechanics and which of the types of these sets can be of positive measure at present this question still remains open). To accumulate specific information we organized a seminar on the study of individual examples. My ideas concerning this topic and closely related problems aroused wide response among young mathematicians in Moscow. [Tikhomirov 1991, p. 521]
\\

\leftskip=0cm
\noindent
The research he was presenting in this important scientific meeting connected to mathematical research in a period in which attention classical mechanics declined in favour of the new mechanics [Dugas 1957]. Thus, his interest could came as a surprise. A contribution that he had just published in the proceedings of the USSR Academy of Sciences was the starting point of research by his students Vladimir Arnold and Yakov G. Sinai in the late 1950s and 1960s. What was the origin of his interest in classical mechanics, and particularly in aspects regarding celestial mechanics? Kolmogorov's views on mathematics, as discussed in successive editions of the entry “Mathematics” in the Soviet Encyclopedia, was tigthly link to 19th century vision of mathematics, that he evoked in the opening words of his Amsterdam lecture:
\\

\leftskip=1cm
\noindent
For the mathematics of the 19th century one of these focal points was the problem of integration of systems of differential equations of classical mechanics, where problems of mechanics, the theory of differential equations, problems of the calculus of variations, multidimensional differential geometry, the theory of analytic functions, and the theory of continuous groups were organically interwoven. After the work of H. Poincaré, the fundamental role of topology for this range of problems became clear. On the other hand, the Poincaré-Carathéodory recurrence theorem initiated the "metrical" theory of dynamical systems in the sense of the study of properties of motions holding for "almost all" initial states of the system. [Kolmgorov 1957/1954, p.355-356].
\\

\leftskip=0cm
\noindent
What circumstances lead Kolmogorov to rethink and elaborate John von Neumann's contributions to classical mechanics – not a central aspect of his extraordinary work – and on work by 1937 Krylov with his collaborator Bogoliuvov on nonlinear mechanics? This is question that Arnold was the first to ask himself. Following his testimony, he both put forward and explanation of his own and ask directly Kolmogorov, obtaining what were a puzzling answer . He has described the conversation twice in published works; in the second one he indicates the date of 1984 for this short exchange between former mentor and former student. In this paper Arnold's account is analyzed in the context of Kolmogorov life and work, from his early years – the final years of the Russian Empire under the Tzarist monarchy – to the years of Stalinist terror, ending with Stalin's death some months before a first paper by Kolmogorov \textit{On Dynamical systems with an integral invariant on the torus}, published on November 31, 1953.

\section{Testimony of a student: A short conversation between Vladimir Igorevich Arnold and Kolmogorov in 1984}

Vladimir Igorevich Arnold was a Ukrainian-born mathematician, born in Odessa on June 12, 1937 and grew up in Moscow\footnote{The Ukrainian origin is due to Arnold's mother, Nina Alexandrova Isakovich, who was of Jewish family. Her mothers family lived in Odessa, where Arnold himself was born in 1937 (the city was then in the USSR Ukrain Soviet Socialist Republic), even if he grew in Moscow.}.
\\
\indent
He entered Moscow State University in 1954 as an undergraduate at the Faculty of Mechanics and Mathematics. He himself will say that he was lucky enough to be the right age at the right time: 
\\

\leftskip=1cm
\noindent
I entered the Faculty for Mechanics and Mathematics of the Moscow State University in 1954 (before Stalin’s death in 1953 or after the invasion to Czechoslovakia in 1968, this would probably have been impossible for me because my mother was a Jew while my grandfather was shot dead in 1938 on the flagrantly false charge of espionage for England, Germany, Greece, and Japan). [Sevryuk 2014, p 3].
\\

\leftskip=0cm
\noindent
The indissoluble union with Kolmogorov manifested itself from the early years: in 1959 he defended his thesis under the supervision of Kolmogorov and in 1961 he received the title of "candidate in physical-mathematical sciences", analogous to the PhD. in the West, at the Keldysh Applied Mathematics Institute in Moscow, with the dissertation containing his famous resolution of Hilbert's 13th problem.
\\
\indent
At the age of just 28, he became a Professor in the Faculty of Mechanics and Mathematics at Moscow State University.

\leftskip=0cm
\noindent
A relationship was established between Arnold and Kolmogorov that goes beyond the pairing of pupil-teacher. In [Arnold 2000] some letters sent to him by Kolmogorov were also published, in which the confidential relationship between the two mathematicians is evident.
\\
\indent
This allows us to consult some direct testimonies which, although sometimes laudatory and not entirely objective, gather some information necessary for the purposes of our research, if analyzed in a critical and rigorous way.
\\
\indent
Kolmogorov published in just ten months, from November 1953 to August 1954, three articles on classical mechanics\footnote{The first two articles [Kolmogorov 1953] and [Kolmogorov 1954] are published in the journal Doklady Akademii Nauk, SSSR, the last one is found in the proceedings of the ICM of 1954, [Kolmogorov 1957/1954]}. In [Kolmogorov 1954] we find the theorem enunciated by the author on the persistence of invariant tori under "small" perturbations, which was the result on which the KAM theory developed in the following years.
\\
\indent
At the time, Kolmogorov was already an established mathematician, having trained with Luzin\footnote{Nikolai Nikolaevič Luzin (1883-1950) was one of the founding fathers of the Moscow school of mathematics. He made important contributions to the foundations of mathematics, measurement theory and topology.} and having distinguished himself internationally with his contributions on probability in the 1930s.
\\
\indent
In September 1954, invited as a speaker in one of the plenary conferences of the International Congress of Mathematicians in Amsterdam, he chose to present the results obtained in the field of classical mechanics, recently published.
These were undoubtedly impactful, as they contradicted a conjecture by Poincaré referred to in the context of celestial mechanics. Although in the first two articles Kolmogorov did not reveal a preponderant interest in celestial mechanics, in the proceedings [Kolmogorov 1957/1954] his interest in this area appears evident, both for a substantial part of the bibliography dedicated to articles from the 1930s on the problem of the three bodies and on celestial mechanics in general and for the numerous references to them in the text.
\\
\indent
In an attempt to trace the cultural origins of Kolmgorov's results on invariant tori and provide a possible answer as to why he was interested in a problem of classical mechanics, so intertwined with questions of celestial mechanics, two excerpts from Arnold are relevant, [Arnold 1997] and [Arnold 200]. In these the author reports similar information, with subtle differences that will be analysed. They were published three years apart, and in any case after Kolmogorov's death in 1987.
\\
\indent
In [Arnold 1997] the author reports, in the form of recollections, the origins of his works on quasi-periodic motions in dynamical systems:
\\

\leftskip=1cm
He later related that he had been thinking about this problem for decades starting from his childhood when he had read Flammarion’s Astronomy, but the success had come only after Stalin’s death in 1953 when a new epoch had begun in the Russian life. The hopes this death raised had a deep impact on Kolmogorov, and the years 1953-1963 were one of the most productive periods in his life.
\\

\leftskip=0cm
\noindent
Instead, the second extract [Arnold 2000] is taken from a volume of the History of Mathematics series, \textit{Kolmogorov in perspective} [AA.VV. 2000], which contains some testimonies of his private life written by students, including Arnold, and Kolmogorov's colleagues.
\\
\indent
Arnold tried to give himself and to give others an answer regarding Kolmogorov's interest - it would seem sudden - towards the study of classical mechanics. 
\\
\indent
Arnold attempted to answer the question himself, but his theory was later invalidated by Kolmogorov himself after a few years. Here is the complete passage in which Arnold's idea and the response of the now friend-colleague are manifested:
\\

\leftskip=1cm
I constructed for myself a theory of the origin of Andrei Nikolaevich's work on invariant tori: it began with his studies of turbulence. In the well-known work of Landau \footnote{He refers to Lev Davidovič Landau (1908 -1968), Soviet physicist winner of the Nobel Prize in physics in 1962, one of the most important physicists of the 20th century. He wrote numerous treatises on mechanics, hydrodynamics, quantum physics and physical statistics.} (1943) it was invariant tori—attractors in the phase space of the Navier-Stokes equation—that were used to "explain" the onset of turbulence.[...] In a discussion at the Landau seminar Andrei Nikolaevich remarked that a transition to an infinite dimensional torus and even to a continuous spectrum can already take place for a finite Reynolds number. On the other hand, even if the dimension of the invariant torus remains finite for a fixed Reynolds number, the spectrum of a conditionally periodic motion on a torus of sufficiently high dimension contains so many frequencies that it is practically indistinguishable from a continuous spectrum. The question as to which of these two cases actually holds was asked more than once by Andrei Nikolaevich. A program for the seminar on the theory of dynamical systems and hydrodynamics was posted on a bulletin board in the Mechanics and Mathematics Department of Moscow State University at the end of the 1950's [...]. Andrei Nikolaevich chuckled about the tori of Landau: "He (Landau) evidently did not know about other dynamical systems." The transition from the tori of Landau to dynamical systems on a torus would be a completely natural train of thought. In the final analysis I almost believed in my theory and (in 1984) asked Andrei Nikolaevich whether it was really so. "No," he answered, "I was not at all thinking of that at the time. The main thing was that there appeared to be hope in 1953. From this I felt an extraordinary enthusiasm.
I had thought for a long time about problems in celestial mechanics, from childhood, from Flammarion, and then—reading Charlier, Birkhoff, the mechanics of Whittaker, the work of Krylov and Bogolyubov, Chazy, Schmidt. I had tried several times, without results. But here was a beginning." [Arnol'd 2000, p 89-90].
\\

\leftskip=0cm
\noindent
Arnold himself probably underestimated how interested Kolmogorov was in celestial mechanics. From the reading of the two testimonies, Kolmogorov's interest in astronomy that dates back to his childhood emerges: he speaks of the famous astronomer and popularizer Camille Flammarion, author of many multi-translated bestsellers. This circumstance underlines the link with astronomy - and therefore with celestial mechanics.
\\
\indent
Other insights emerge: in both excerpts he speaks of a hope reborn since 1953. In particular, in the first it emerges more directly that it is connected with Stalin's death.
\\
\indent
Other insights emerge. In both excerpts he speaks of a hope reborn since 1953. In particular, in the first it emerges more directly that it is connected with Stalin's death.
In the period 1936 - March 1953 the population of the Soviet Union was bent by Stalin's Great Purges and the terror policy implemented. Internal tensions, the growing threat of a Second World War, the Iron Curtain, helped create distances and barriers between the USSR and the rest of the world.
\\
\indent
In the mathematical academic field, starting from the second half of the 1930s, connections with foreign countries were interrupted. Travel and attendance at conferences were banned\footnote{In Deminov [2009], the author elaborates on abrupt interruptions of travel and relations with France.}. 
\\
\indent
Furthermore, in the years 1936-37 there was a real witch hunt in the world of astronomy, which involved not only astronomers but various scientific sectors, such as astronomers, geologists, geophysicists, geodesists. The purge devastated the powerful group of astronomers (more than two dozen) active in a major network of observatories that developed in Russia's last period under the Tsar\footnote{See [McCutcheon 1991], [Eremeeva 1995]}. This dramatic event may have prompted Kolmogorov to wait until 1953 to devote himself to - or simply publish - articles on classical mechanics linked to Poincaré's works on celestial mechanics.

\section{Kolmogorov and celestial mechanics}

In [Kojevnikov 2008, p 115-116], the author discusses the beginning of the great Soviet experiment which, starting from 1860, led science in a few decades to such a development as to be competitive with the most advanced European countries such as Germany and Great Britain: 
\\

\leftskip=1cm
To begin with, the Soviets mounted their belief in science on top of a preexisting and rather high foundation. The cult of science flourished across Europe at the beginning of the twentieth century. It happened to be particularly prominent in the Russian empire, which had only recently embarked upon industrialization and modernization.
Almost all parts of the political spectrum bought into it, although for different reasons. For Russian liberals, science was synonymous with economic and social progress; for the radical intelligentsia, including the yet utterly insignificant and marginal Bolsheviks on the very left, it was the closest ally of the revolution. Many among the monarchists, too, placed high hopes on modern science as a remedy for the country’s relative economic backwardness vis-à-vis Germany, France, and Britain (other European countries rarely figured in the comparison). After the Great Reforms of the 1860s, they helped institutionalize science and promote the research imperative at Russian universities, hoping that at the very least it could distract unruly students from pursuing dangerous political temptations. 
\\

\leftskip=0cm
\noindent
Kolmogorov lived his childhood in the midst of the revolution and there is no doubt that he absorbed all the influences of those years.
Others to numerous biographical and autobiographical notes that testify to it\footnote{Biographical essays, [Shiryaev 1992] and [Shiryaev 2000] draw many informations from Kolmogorov's reminiscences included in his posthumous book on mathematics as science and as a profession [Kolmogorov 1988]}, his initial indecision between mathematics and metallurgy also confirms it. In 1963, in an interview with Ogonek, one of Russia's oldest illustrated weeklies, founded in 1899, he recalled that "engineering was then perceived as something more serious and necessary than pure science." This thought is clearly the son of the historical and cultural period lived in Russia.
\\
\indent
As we have read from Arnold, his propensity towards applied sciences also manifested itself in astronomy.
The references in the previous section take us back to the astronomer Camille Flammarion (1842–1925), a famous French astronomer, publisher and science popularizer. A prolific and multi-translated author, during his career he published more than fifty works, among which the most famous were popular astronomy guides. To get an idea of the impact it had in France, and in the world in general, let's read a few lines reported in an obituary written by the English astronomer William Porthoue (1877-1964), a member of the Manchester Astronomical Society from 1905 until his death and editor from 1913 to 1924 of the Journal of the Manchester Astronomical Society:
\\

\leftskip=1cm
Camille Flammarion might be described as the apostle of popular astronomy. His numerous literary works had for object primarily the popularisation of astronomical study in all its manifolds branches [...]. 
Flammarion was not content to spread abroad the gospel of astronomy by book and pamphlet. He believed in the practical application of his theories for the spread of a universal knowledge of the sky. [Porthouse 1925, p 951]
\\

\leftskip=0cm
\noindent
Strongly convinced that the study of science was for everyone, Flammarion collaborated with a large number of magazines and newspapers, actively participating in the great scientific emancipation movement of the second half of the nineteenth century. His books are rich in figures and illustrations and are written with direct and persuasive communication, with a style capable of enthralling and enthusing the reader.
But which \textit{Flammarion's Astronomy} could Kolmogorov have read? 
\\
\indent 
Analyzing the time span in which he as a child would have read the works and translations in Russian, we could restrict the field to two possibilities:  one of which is the famous \textit{Astronomie Populaire}, [Flammarion 1880] published in 1880 in Paris by the publishing house "C. Marpon et E. Flammarion" \footnote{Ernest Flammarion (1846 - 1936), French publisher and Camille's brother.}, translated for the first time into Russian as early as 1897 and in various subsequent editions. Divided into six chapters - The earth, The moon, The sun, The planetary worlds, The comets, The stars - it is intended to be a book aimed at everyone to teach the elementary knowledge of astronomy in an extremely didactic and popular form.
The opening words of chapter one is an introduction to the entire book:
\\

\leftskip=1cm
Ce livre est écrit pour tous ceux qui aiment a se rendre compte des choses qui les entourent, et qui seraient heureux d'acquérir sans fatig-ue une notion élémentaire et exacte de l'état de l'univers. 
N'est-il pas agréable d'exercer notre esprit dans la contemplation des grands spectacles de la nature? N'est-il pas utile de savoir au moins sur quoi nous marchons, quelle place nous occupons dans l'infini, quel est ce soleil dont les rayons bienfaisants entretiennent la vie terrestre, quel est ce ciel qui nous environne, quelles sont ces nombreuses étoiles qui pendant la nuit obscure répandent dans l'espace leur silencieuse lumière? Cette connaissance élémentaire de l'univers, sans laquelle nous végéterions comme les plantes, dans l'ignorance et l'indifférence des causes dont nous subissons perpétuellement les effets, nous pouvons l'acquérir, non-seulement sans peine, mais encore avec un plaisir toujours grandissant. Loin d'être une science isolée et inaccessible, l'Astronomie est la science qui nous touche de plus près, celle qui est la plus nécessaire à notre instruction générale, et en même temps celle dont l'étude offre le plus de charmes et garde en réserve les plus profondes jouissances.\footnote{This work is written for those who wish to hear an account of the things which surround them, and who would like to acquire, without hard work, an elementary and exact idea of the present condition of the universe.
It is not pleasant to exercise our minds in the contemplation of the great spectacles of nature? It is not useful to know, at least, upon what we tread, what place we occupy in the infinite, the nature of the sun whose rays maintain terrestrial life, of the sky which surrounds us, of the numerous stars which in the darkness of night scatter through space their silent light? This elementary knowledge of the universe, without which we live, like plants, in ignorance and indifference to the causes of which we perpetually witness the effects, we can acquire not only without difficulty, but with an ever-increasing pleasure. Far from being a difficult and inaccessible science, Astronomy is the science which concerns us most, the one most necessary for our general instruction, and at the same time the one which offers for our study the greatest charm and keeps in reserve the highest enjoyments. 
\\
Populair Astronomy (1894), English version, translated by J. Ellard Gore, London, Chatto $\&$ Windus, Piccadilly, p.1}
\\

\leftskip=0cm
\noindent
Another probable reading, although less famous than the first but well known in the innovative circles of that time, is \textit{Initiation Astronomique} [Flammarion 1908]. It was a work aimed at children, also with numerous illustrations, published in 1908 in Paris by Libraire Hachette et C$^{Ie}$ and translated into Russian in the same year, when Kolmogorov was five years old and was in Tunoshna where he was studying in the innovative school of aunts.
\\
The booklet was published in the series of "Initiations scientifiques" directed by the mathematician Charles-Ange Laisant \footnote{Charles-Ange Laisant (1841-1920). French mathematician and politician, he was a deputy from Nantes and professor at the École polytechnique in Paris. He dealt with mechanics, geometry and algebra and, mainly, with the teaching of mathematics and related reform.} and, as he himself writes in the opening introductory pages, "Il est destiné, entre le mains de l'éducateur, à servir de guide pour la formation de esprit des tout jeunes enfants - de quatre à douze ans - afin de meubler leur intelligence de notions saines et justes, et de les préparer ainsi à l'étude, qui viendra plus tard." \footnote{Eng. tr.: It is intended, in the hands of the educator, to act as a guide for the formation of the mind of very young children - from four to twelve years old - in order to provide their intelligence with sound and correct notions, and thus prepare them for examination, which will come afterwards.} [Flammarion 1908, p. V].
\\
We also find a brief introduction by the author who, in his words, expresses all his passion in this project and affirms the centrality of astronomy in scientific thought: 
\\

\leftskip=1cm
J'ai toujours pensé aussi qu'il n'est pas nécessaire d'ennuyer le lecteur puor l'instruire, et que si pendant tant de siècles, l'Astronomie, la plus belle des sciences, celle qui nous apprend où nous sommes et qui nous dévoile les splendeurs de l'Univers, est restée à peu orès ignorée de l'immense majorité des habitants de notre  planète, c'est parce qu'elle a toujours été mal enseignée dan les Ècoles. Aujourd'hui, enfin, on commence à la trouver intéressante, à lire le grand livre de la Nature, à vivre un peu plus intellectuellement. \footnote{Eng. tr.: I have also always thought that it is not necessary to bore the reader to instruct him, and that if for so many centuries astronomy, the most beautiful of the sciences, the one that teaches us where we are and that reveals the splendors of the Universe, has remained almost ignored by the vast majority of the inhabitants of our planet, it is because it has always been badly taught in schools. Today, finally, we begin to find it interesting, to read the great book of Nature, to live a little more intellectually.} [Flammarion 1908, p. VII]
\\

\leftskip=0cm
\noindent
Although it is not possible to ascertain with certainty which of the cited texts Kolmogorov read, it is clear that this author played a fundamental role in the birth of the interest of the child Kolmogorov in the stars and celestial mechanics.

\section{Final remarks}

The title of the plenary lecture given by Kolmogorov at the ICM in 1954, \textit{The general theory of dynamical systems and classical mechanics}, referred to both the centennial tradition of the mathematical study of motion, elasticity and a growing number of physical phenomena, than to the theory of dynamical systems of the twentieth century, rooted in classical mechanics but at the same time a paradigmatic example of new mathematical approaches developed after 1900.
\\
\indent
The adjective “classical” served to mark the separation between the tradition of “rational mechanics” – from its Newtonian source to the reformulations of Joseph Louis Lagrange (1788) and William Rowan Hamilton (1833) – and the new quantum mechanics.
\\
\indent
The theoretical and epistemological crisis of classical mechanics had distanced it from the center of the mathematical scene in the first half of the twentieth century, also due to the vitality of modern algebra and of branches of mathematics operating in abstract universes without any connection with physical phenomena, developments technology and other applications.
\\
\indent
Kolmogorov's clear intention was to bring back to this center "the problem of integrating the systems of differential equations of classical mechanics", in his own words a "focal point for nineteenth-century mathematics". During the conference he discussed the proof of a theorem which he offered as a contribution in that direction. The theorem he proved and presented in Amsterdam was published a few days before the plenary conference, on August 31, in the USSR mathematical journal \textit{Doklady Akademi Nauk}, in an article entitled \textit{On the conservation of conditionally periodic motions under small variations of the Hamilton function}. This theorem, or rather the research program it contained - as he demonstrated to the international mathematical community during his 1954 conference - was the starting point of the KAM theory, of great importance for the history of mathematics and science, because, as recently stated by Scott Dumas in [Dumas 2014], the KAM theory has made it possible to complete the “true picture of classical mechanics – often thought to have been sketched in the 17th century”:
\\

\leftskip=1cm
\noindent
Right from the start, after enunciating his laws of mechanics and gravitation, Isaac Newton ran into difficulties using those laws to describe the motion of three bodies moving under mutual gravitational attraction (the so-called ‘three body problem’). For the next two centuries, these difficulties resisted solution, as the best minds in mathematics and physics concentrated on solving other, increasingly complex model systems in classical mechanics (in the abstract mathematical setting, to ‘solve’ a system means showing that its trajectories move linearly on so-called ‘invariant tori’). But toward the end of the 19th century, using his own new methods, Henri Poincaré confronted Newton’s difficulties head-on and discovered an astonishing form of ‘unsolvability,’ or chaos, at the heart of the three body problem. This in turn led to a paradox. According to Poincaré and his followers, most classical systems should be chaotic; yet observers and experimentalists did not see this in nature, and mathematicians working with model systems could not (quite) prove it to be true either. The paradox persisted for more than a half century, until Andrey Kolmogorov unraveled it by announcing that, against all expectation, many of the invariant tori from solvable systems remain intact in chaotic systems. These tori make most systems into hybrids – they are a strange, fractal mixture of regularity and chaos.
\\

\leftskip=0cm
\noindent
KAM theory - says Dumas - is a breakthrough, a revolution, comparable to quantum mechanics and the theory of relativity, although it still lacks general understanding today.
\\
\indent
Kolmogorov dealt with conservative dynamical systems, which formed an important part of 19th-century classical mechanics, in particular perturbed Hamiltonian systems with "small perturbations", a central topic in celestial mechanics - think of the three-body problem. In the wake of Poincaré's studies and his \textit{Theorem of non-existence of motion} of 1887, but above all after the introduction of ergodic theory by Boltzmann, mathematicians and physicists were practically convinced that perturbed Hamiltonian systems were all chaotic, even ergodic.
Kolmogorov's theorem shows that this is false for the considered perturbed Hamiltonian systems, under certain conditions.
\\
\indent
Addressing the international mathematical community, the Russian mathematician wanted to show an exceptional example of unexpected and profound relationships between different branches of mathematics: this was an openly stated reason for his choice. The keynote address of his Amsterdam conference clearly suggests that this theorem was the core of an entire research program, which he felt the urgency to share after presenting the mathematical core of his findings in the two Doklady papers.

\newpage

\section*{Bibliography}
\addcontentsline{toc}{section}{Bibliography}
[AA.VV. 2000] \textsc{AA. VV.} (2000) \emph{Kolmogorov in perspective}. Providence, R.I. American Mathematical Society, London Mathematical Society.
\\[1em]
[Abramov 2010] \textsc{Abramov} Alexander, (2010) \emph{Toward a History of Mathematics Education Reform in Soviet Schools
(1960s–1980s)} in [Karp, Vogeli (2010)], pp. 87-140. 
\\[1em]
[Arnold 1997] \textsc{Arnol'd} Vladimir Igorevič (1997) \emph{Vladimir Igorevič Arnold. Selected-60}. Phasis, Moscow, pp. 727-740. (in Russian). 
[Eng. tr. in  \textsc{Sevryuk M.B.} (2014) \emph{From Superpositions to KAM Theory}, Regular and Chaotic Dynamics 19(6) pp. 734-744].
\\[1em]
[Arnold 2000] \textsc{Arnol'd} Vladimir Igorevič (2000) \emph{On A. N. Kolmogorov}. [in AA.VV. \emph{Kolmogorov in perspective}], pp. 89-108.
\\[1em]
[Arnold 2004] \textsc{Arnol'd} Vladimir Igorevič (2004) \emph{From Hilbert's Superposition Problem to Dynamical Systems}. The American Mathematical Monthly 111(7), pp. 608-624.
\\[1em]
[Arnold 2007] \textsc{Arnol'd} Vladimir Igorevič (2007) \emph{Yesterday and long ago}. Springer/Phasis.
\\[1em]
[Charpentier, Lesne, Nikolski 2004] \textsc{Charpentier} Eric, \textsc{Lesne} Annick, \textsc{Nikolski}  Nikolai Kapitonovich (eds.) (2004) \emph{L'héritage de Kolmogorov en mathématiques}. Paris, Éditions Belin [eng. tr. Berlin Heidelberg, Springer Verlag 2007].
\\[1em]
[Deminov 2009] \textsc{Deminov} Sergej Sergeevich (2009) \emph{Les relations mathématiques Franco-Russes entre les deux guerres mondiales}. Revue d'histoire des sciences, 2009/1 Tome 62, pp. 119-142.
\\[1em]
[Diacu, Holmes 1996] \textsc{Diacu} Florin \textsc{Holmes} Philip (1996) \emph{Celestial encounters: The origin of chaos and stability}. Princeton, New York, Pricenton University Press.
\\[1em]
[Diner 1992] \textsc{Diner} Simon (1992) \emph{Les voies du chaos déterministe dans l'école russe}. In [Dahan Dalmedico, Chabert, Chemla 1992], pp 331-370.
\\[1em]
[Dugas 1957] \textsc{Dugas} René (1957) \emph{A history of mechanics}. London, Routledge $\&$ Kegan Paul LTD. [Trasleted into English by J.R.Maddox of \textsc{Dugas} René (1950)  \emph{Histoire de la mécanique}, Switzerland, Neuchâtel, Éditions du Griffon]
\\[1em]
[Dumas 2014] \textsc{Dumas} H. Scott (2014) \emph{The KAM story a friendly introduction to the content, history, and significance of classical Kolmogorov-Arnold-Moser theory}.
Singapore, World Scientific Publishing.
\\[1em]
[Eremeeva 1995] \textsc{Eremeeva} A. I. (1995) \emph{Political repression and personality: the hostory of political repression against soviet astronomers}. Journal for the History of Astronomy 26(4), pp 297-324.
\\[1em]
[Flammarion 1880] \textsc{Flammarion} Camille (1880) \emph{Astronomie Populaire: Description générale du ciel}. Paris, C. Marpon et E. Flammarion Éditeurs.
\\[1em]
[Flammarion 1908] \textsc{Flammarion} Camille (1908) \emph{Initiation Astronomique}. Collection Des Initiations Scientifiques fondée par C.-A. Laisant, Paris, Librairie Hachette et C$^Ie$.
\\[1em]
[Gerretsen, De Groot 1957] \textsc{Gerretsen} Johan C.H., \textsc{De Groot} Johannes (eds.) (1957) \emph{Proceedings of the International Congress of Mathematicians 1954 (Amsterdam September 2 - 9)}. Groningen/Amsterdam, 
Erven P. Noordhoff N.V./North-Holland Publishing Co.
\\[1em]
[Gordin, Hall, Kojevnikov 2008] \textsc{Gordin} Michael, \textsc{Hall} Karl, \textsc{Kojevnikov} Alexei (eds.) (2008) \emph{Intelligentsia Science: The Russian Century, 1860-1960}. Chicago, Chicago University Press, Osiris, Vol.23(1).
\\[1em]
[Graham 1993] \textsc{Graham} Loren R. (1998) \emph{Science in Russia and the Soviet Union: A Short History}. New York, Cambridge University Press.
\\[1em]
[Graves, Hille, Smith, Zariski 1955] \textsc{Graves} Lawrence M., \textsc{Hille} Einar, \textsc{Smith} Paul A., \textsc{Zariski} Oscar (eds.) (1955) \emph{Proceedings of the International Congress of Mathematicians 1950 (Cambridge, Massachusetts, U.S.A. 1950}. Providence, RI, American Mathematical Society. 
\\[1em]
[Karp 2012] \textsc{Karp} Alexander (2012) \emph{Soviet mathematics education between 1918 and 1931: a time of radical reforms}. Karlsruhe, ZDM Mathematics Education, vol. 44, pp. 551–561.
\\[1em]
[Karp 2014] \textsc{Karp} Alexander, (2014) \emph{Mathematics Education in Russia} in [Karp, Schubring 2014] pp. 303-322.
\\[1em]
[Karp, Schubring 2014] \textsc{Karp} Alexander, \textsc{Schubring} Gert (edited by) (2014) \emph{Handbook on the History of Mathematics Education}. New York, Springer.
\\[1em]
[Karp, Vogeli 2010] \textsc{Karp} Alexander,\textsc{Vogeli} Bruce (eds.) (2010) \emph{Russian mathematics education: history and world significance}. Series of mathematics educations, vol. 4 London/New Jersey/Singapore, World Scientific Publishing Co Pte Ltd.
\\[1em]
[Kojevnikov 2002] \textsc{Kojevnikov} Alexei (2002) \emph{Introduction: A New History of Russian Science}. Cambridge University Press. Science in Context 15(2), pp.177–182.
\\[1em]
[Kojevnikov 2008] \textsc{Kojevnikov} Alexei (2008) \emph{The Phenomenon of Soviet Science}. In [Gordin, Hall, Kojevnikov 2008], pp. 115-135.
\\[1em]
[Kolmogorov 1923] \textsc{Kolmogorov} Andrei Nikolaevich (1923) \emph{Une série de Fourier-Lebesgue divergente presque partout}. Warsaw, Fundamenta Mathematicae 4, pp. 324-328.
\\[1em]
[Kolmogorov 1953] \textsc{Kolmogorov} Andrei Nikolaevich (1953) \emph{On Dynamical systems with an integral invariant on the torus} (in Russian). Doklady Akademii Nauk SSSR, 93(5), pp. 763-766. [English translation in: \emph{Selected works of A.N. Kolmogorov}, vol.I edited by V. M. Tikhomirov Dordrecht, Kluwer Academic Publishers, 1991, pp. 344-348].
\\[1em]
[Kolmogorov 1954] \textsc{Kolmogorov} Andrei Nikolaevich (1954) \emph{On the preservation of conditionally periodic motions under small variations of the Hamilton function } (in Russian). Doklady Akademii Nauk SSSR, 98(4), pp. 527-530.
[English translation in: \emph{Selected works of A.N. Kolmogorov}, vol.I edited by V. M. Tikhomirov Dordrecht, Kluwer Academic Publishers, 1991, pp. 349-354].
\\[1em]
[Kolmogorov 1957/1954] \textsc{Kolmogorov} Andrei Nikolaevich (1957/1954) \emph{The general theory of dynamical systems and classical mechanics}, in [Tikhomirov 1991, vol.1, pp.355-374]. Original edition in Russian in \textsc{Gerretsen Johan C.H., De Groot Johannes} (1957) \emph{Proceedings of the International Congress of Mathematicians 1954 (Amsterdam September 2 - 9)}, North Holland, Amsterdam, vol. 1, pp. 315–333.
\\[1em]
[Krylov,Bogolioubov 1937] \textsc{Krylov} Nikolai Mitrofanovich, \textsc{Bogolioubov} Nikolaj Nikolaevich (1937) \emph{La Theorie Generale De La Mesure Dans Son Application A L'Etude Des Systemes Dynamiques De la Mecanique Non Lineaire}. Princeton University, Annals of Mathematics, Second Series, vol. 38(1), pp. 65-113.
\\[1em]
[Kutateladze 2012] \textsc{Kutateladze} Semën Samsonovich (2013) \emph{The Tragedy of Mathematics in Russia.}. Siberian Electronic Mathematical Reports, Vol. 9, pp. A85–A100.
\\[1em]
[Kutateladze 2013] \textsc{Kutateladze} Semën Samsonovich (2013) \emph{An epilog to the Luzin case}. Siberian Electronic Mathematical Reports, Vol. 10, pp. A1–A6.
\\[1em]
[Levin 1990] \textsc{Levin}  Aleksey E. (1990) \emph{Anatomy of a Public Campaign: "Academician Luzin's Case" in Soviet Political History}. Slavic Review Vol. 49(1), pp. 90-108
\\[1em]
[Mazliak 2018] \textsc{Mazliak} Laurent (2018) \emph{The beginnings of the Soviet encyclopedia. The utopia and misery of mathematics in the political turmoil of the 1920s}. Centaurus 60(1-2), pp 25-51.
\\[1em]
[McCutcheon 1991] \textsc{McCutcheon} Robert A. (1991) \emph{The 1936-1937 Purge of Soviet Astronomers}. Cambridge University Press, Spring 50(1), pp. 100-117. 
\\[1em]
[Porthouse 1925] \textsc{Porthouse} William, (1925) \emph{Obituary: M. Camille Flammarion}. Nature 115(2903), pp.951-952.
\\[1em]
[Shiryaev 1989] \textsc{Shiryaev} Albert Nikolayevich (1989) \emph{Kolmogorov: Life and Creative Activities}. The Annals of probability 17(3) pp. 866-944.
\\[1em]
[Shiryaev 1992] \textsc{Shiryaev} Albert Nikolayevich (1992) (ed.) \emph{Selected Works of A. N. Kolmogorov, vol II}. Dordrecht, Springer Science+Business Media.
\\[1em]
[Shiryaev 2000] \textsc{Shiryaev} Albert Nikolayevich (2000) \emph{Andrei Nikolaevich Kolmogorov (April 25, 1903 to October 20, 1987) A Biographical Sketch of His Life and Creative Paths}. in [AA.VV. 2000], pp. 1-88.
\\[1em]
[Tikhomirov 1991] \textsc{Tikhomirov} Vladimir Mikhailovich (1991) \emph{Selected works of A.N. Kolmogorov, vol. I}. Dordrecht, Kluwer Academic Publishers.
\\[1em]
[Zdravkovska, Duren 2007] \textsc{Zdravkovska} Smilka, \textsc{Duren} Peter Larkin (eds.) (2007) \emph{Golden Years of Moscow Mathematics}.  Second edition, American Mathematical Society, History of Mathematics, vol.6.

\end{document}